\documentclass{amsart}

\usepackage{amsmath,amssymb,amscd,amsfonts}

\newtheorem{theorem}{Theorem}
\newcommand{\bt}{\begin{theorem}}
\newcommand{\et}{\end{theorem}}
\newtheorem{lemma}{Lemma}
\newcommand{\bl}{\begin{lemma}}
\newcommand{\el}{\end{lemma}}
\newtheorem{corollary}{Corollary}
\newcommand{\bc}{\begin{corollary}}
\newcommand{\ec}{\end{corollary}}
\newcommand{\beq}{\begin{equation}}
\newcommand{\eeq}{\end{equation}}
\newcommand{\benum}{\begin{enumerate}}
\newcommand{\eenum}{\end{enumerate}}
\newcommand{\N}{\ensuremath{ \mathbf N }}
\newcommand{\Z}{\ensuremath{\mathbf Z}}

\newcommand{\mca}{\ensuremath{ \mathcal A}}
\newcommand{\mcb}{\ensuremath{ \mathcal B}}
\newcommand{\mcc}{\ensuremath{ \mathcal C}}

\newcommand{\mcm}{\ensuremath{ \mathcal M}}

\title{Additive systems and a theorem of de Bruijn}

\author{Melvyn B. Nathanson}
\address{Lehman College (CUNY),Bronx, New York 10468}
\email{melvyn.nathanson@lehman.cuny.edu}

\subjclass[2010]{Primary 11A05, 11B75.} 
\thanks{Supported in part by a grant from the PSC-CUNY Research Award Program.}

\dedicatory{In memoriam Nicolaas Govert  de Bruijn}

\begin{document}

\maketitle

\begin{abstract}
This paper proves a theorem of de Bruijn that classifies 
additive systems for the nonnegative integers, that is, 
families $\mca = (A_i)_{i\in I}$ of sets of nonnegative integers, 
each set containing 0, such that every nonnegative integer can be written 
uniquely in the form $\sum_{i\in I} a_i$ with $a_i \in A_i$ for all $i$ 
and $a_i \neq 0$ for only finitely many $i$.  
\end{abstract}

\section{Additive systems}

Let \N, $\N_0$,  and \Z\ denote the sets of positive integers, nonnegative integers, 
and all integers, respectively.  For integers $a$ and $b$ with $a < b$, we define the  \emph{intervals of integers} 
$[a,b] = \{ n \in \Z : a \leq n \leq b \}$ and 
$[a,b) = \{ n \in \Z : a \leq n < b \}.
$
For $A \subseteq \Z$ and $g\in \Z$, the \emph{dilation}
\index{dilation} of the set $A$ by $g$ is the set $g\ast A = \{ ga : a \in A\}$.  

Let $I$ be a nonempty finite or infinite set, 
and let $\mca = (A_i)_{i\in I}$ be a family of  sets of integers with $0 \in A_i$ and 
$|A_i| \geq 2$ for all $i \in I$.   
We may also call \mca\ a \emph{sequence} if $I= \N$ or if $I$ is an interval of integers.  
Each set $A_i$ can be finite or infinite.  
We say that a set $X$ \emph{belongs to \mca} if $X = A_i$ for some $i \in I$.
The \emph{sumset} $S = \sum_{i\in I} A_i$ is the set of all integers $n$ 
that can be represented in the form $n = \sum_{i\in I} a_i$, 
where $a_i \in A_i$ for all $i \in I$ and $a_i  \neq 0$ for only finitely many $i \in I$.  
If every element of  $S$ has a \emph{unique} 
representation in the form $n = \sum_{i\in I} a_i$, 
then we call \mca\ a  \emph{unique representation system for $S$}, 
\index{unique representation system}
and we write 
$
S =  \bigoplus_{i\in I} A_i.
$
If \mca\ is a unique representation system for $S$, 
then $A_i \cap A_j = \{ 0 \}$ for all $i \neq j$. 
The condition $|A_i| \geq 2$ for all $i\in I$ implies that  
$A_i = S$ for some $i\in I$ if and only if $|I|=1$.  
Moreover, if $I^{\flat} \subseteq I$ and $S =  \sum_{i\in I^{\flat}} A_i$, 
then $S =  \bigoplus_{i\in I^{\flat}} A_i$ and $I = I^{\flat}$.

The family $\mca = (A_i)_{i\in I}$ is  an \emph{additive system}
\index{additive system} if \mca\  is a unique representation system 
for the set of nonnegative integers.  
Equivalently, \mca\  is an additive system 
if $\N_0 =  \bigoplus_{i\in I} A_i$.  

The object of this paper is to prove a beautiful theorem of deBruijn 
in additive number theory that completely classifies additive systems.   
The  only number theory used in the proof is the \emph{division algorithm}:
For every positive integer $g$ and for every integer $n$ 
there exist unique integers $x$ and $r$ 
with $r \in [0,g)$ such that $n = gx+r$.

\textbf{Example 1:}  
For $g \geq 2$, let
\[
A_1 = \{0,1,2,\ldots, g-1\} = [0,g)
\]
and
\[
A_2 = \{0,g,2g,3g,4g,\ldots\} = g\ast \N_0.
\]
The division algorithm implies that  $\mca = (A_i)_{i\in [1,2]}$ is an additive system.  
More generally, let $\mca = (A_i )_{i \in I}$ be an additive system.  
Let $I_1 = I \cup \{i_1\}$, where $i_1 \notin I$, and define the sets 
$A'_{i_1} = [0,g)$ and $A'_i = g \ast A_i$ for all $i \in I$.  
Again, the division algorithm implies that $\mca' = (A'_i )_{i \in I_1}$ is an additive system. 
We call $\mca'$ the \emph{dilation of the additive system \mca\ by 
the integer $g$}, and we write $\mca' = g\ast \mca$.

\textbf{Example 2:}
For $i=1,2,3,\ldots,$ let 
\[
B_i = \{ 0,2^{i-1} \}  = 2^{i-1}\ast [0,2).
\]
Because every nonnegative integer can be written uniquely 
as a finite sum of pairwise distinct powers of 2, 
the family $\mcb = (B_i )_{i \in \N}$ is an additive system, 
called the \emph{binary number system}.  
More generally, for any integer $g \geq 2$, let 
\[
C_i  = g^{i-1}\ast [0,g)
\]
for $i=1,2,3,\ldots.$
The additive system $\mcc = \{C_i\}_{i \in \N}$ 
is the \emph{$g$-adic number system}.

\textbf{Example 3:}
Let
\begin{align*}
M_1 & = \{0,1,2, 3,\ldots, 11\} = [0,12) \\
M_2 & = \{0,12,24,36,\ldots, 228\} = 12 \ast [0,20)  \\
M_3 & = \{ 0,  240,  480,   720, 960,  \ldots  \}  = 240\ast \N_0.
\end{align*}
Applying the division algorithm  
with $r=2$, $g_1 = 12$ and $g_2 = 20$,  we see that 
$\mcm = (M_i)_{i\in [1,3]}$ is an additive system.  For example, 
\[
835 = 7 + 108 + 720 = 1\cdot 7 + 12\cdot 9 + 240\cdot 3 \in \sum_{i \in [1,3]} M_i.
\]
In pre-1971 British currency, there were 20 shillings in a pound and 12 pence (or pennies) in a shilling, hence 240 pence in a pound.  
Thus, 835 pence were equal to 3 pounds, 9 shillings, and 7 pence.  
The additive system \mcm\ is   \index{British monetary system}  
the old \emph{British monetary system}.

The following result generalizes Example 3.

\bt     \label{deBruijn:theorem:DivisionAlgorithm}
Let $r \in \N$ and let $(g_i)_{i \in [1,r]}$ be a finite sequence
of not necessarily distinct integers such that $g_i \geq 2$ for all $i \in [1,r]$.  
Let $G_0 = 1$ and $G_i = \prod_{j=1}^i g_j$ for $i \in [1,r]$. 
Then
\beq           \label{deBruijn:DA1}
[0, G_r) = \bigoplus_{i \in [1,r] } G_{i-1}\ast [0,g_i)
\eeq
and
\beq           \label{deBruijn:DA2}
\N_0 = \bigoplus_{i \in [1,r]} G_{i-1}\ast [0,g_i) \oplus G_r\ast \N_0.
\eeq
\et

Thus, the family $(G_{i-1} \ast [0,g_i))_{i\in [1,r]}$ 
is a unique representation system for the interval $[0, G_r)$, 
and this family together with the set $G_r\ast \N_0$ is an additive system.

\begin{proof}
The proof is by induction on $r$. 
The case $r=1$ is Example 1.

Let $r \geq 2$ and assume the Theorem holds for $r-1$.
For  $n \in \N_0$ there are unique integers $x_1,\ldots, x_{r-1},x'_r$ with 
$x_i \in [0, g_i)$ for $i \in [1,r-1)$ and $x'_r \in \N_0$ such that 
\[
n = \sum_{i=1}^{r-1} G_{i-1} x_i + G_{r-1}x'_r.
\]
Applying the division algorithm to $x'_r$, we obtain unique integers 
$x_r \in [0,g_r)$ and $x_{r+1} \in \N_0$ such that $x'_r = x_r + g_rx_{r+1}$,  and so  
\begin{align*}
n & = \sum_{i=1}^{r-1} G_{i-1} x_i + G_{r-1} ( x_r + g_r x_{r+1})   = \sum_{i=1}^{r} G_{i-1} x_i + G_r x_{r+1}.
\end{align*}
The inequality
\[
0 \leq \sum_{i=1}^r G_{i-1} x_i \leq \sum_{i=1}^r G_{i-1} (g_i-1)
=  \sum_{i=1}^r G_i -  \sum_{i=1}^rG_{i-1} = G_r-1
\]
implies that $n \in [0,G_r)$ if and only if $x_{r+1} = 0$.  
This completes the proof.
\end{proof}

\section{Dilation and contraction}
In this section we describe two operations on additive systems 
that produce new additive systems.
Let $\mca = (A_i)_{i\in I}$ be an additive system,
Without loss of generality, and for simplicity of notation, 
we shall assume that $I \cap \N = \emptyset$.

In Example 1 we described the dilation 
of the additive system by an integer $g \geq 2$.
We define \emph{dilation by a finite family}
$(g_i)_{i\in [1,r]}$ of integers $g_i \geq 2$ by iterated dilation by integers:   
\[
(g_i)_{i\in [1,r] } \ast \mca 
 = g_1 \ast \left(   g_2 \ast \left( \cdots  \ast \left(
g_{r-1} \ast \left(    g_r \ast \mca    \right)  \right)
\cdots   \right)    \right)  = (A'_i)_{i\in [1,r]\cup I}
\]
where 
\[
A'_i = 
\begin{cases}
G_{i-1}\ast [0,g_i) & \text{ if } i\in [1,r]\\
G_r \ast A_i & \text{ if } i \in I.
\end{cases}
\]
and $G_0 =1$ and $G_i = \prod_{j\in [1,i]} g_j$ for $i \in [1,r]$. 

Note that dilation of additive systems by finite families of integers is 
not commutative.  
For example, if $g_1 \neq g_2$, then 
$g_1\ast ( g_2 \ast \mca)$ consists of 
$\left( g_1 g_2\ast A_i \right)_{i\in I}$ and the sets 
$[0,g_1)$ and $g_1\ast [0,g_2)$, while 
$g_2 \ast (g_1 \ast \mca)$ consists of the sets 
$\left( g_1 g_2\ast A_i \right)_{i\in I}$ and the sets 
$[0,g_2)$ and $g_2\ast [0,g_1)$.  
Because $[0,g_1) \neq [0,g_2)$, it follows that 
$(g_i)_{i\in [1,2] } \ast \mca \neq (g_{3-i})_{i\in [1,2] } \ast \mca$.

The following two lemmas  follow immediately from the definition of dilation 
and the definition of additive system, respectively.  
The first lemma shows that the dilation of a dilation is a  dilation, 
or, equivalently, that dilation is associative. 
The second shows that partitioning an index set  produces a new  additive system.

\bl     \label{deBruijn:lemma:DilationTransitivity}
Let $\mca$, $\mcb$, and $\mcc$ be additive systems.  
If the additive system \mca\ is a dilation of the additive system \mcb\ 
 by the finite sequence $(g_i)_{i\in [1,r]}$,
and if \mcb\ is a dilation of the additive system \mcc\
by the finite sequence $(g'_j)_{j \in [1,s]}$,
then \mca\ is a dilation of the additive system \mcc\ dilated by $(g_i)_{i\in [1,r+s]}$,
where $g_{r+j} = g'_j$ for $j \in [1,s]$.
\el

\bl     \label{deBruijn:lemma:contraction}
Let $\mcb = (B_j )_{j\in J} $  be an additive system.  
If $\{J_i\}_{i\in I}$ is a partition of $J$ into pairwise disjoint nonempty sets, and  if 
\[
A_i= \sum_{j \in J_i} B_j
\]
then $\mca = (A_i )_{i\in I}$ is an additive system.
\el

An additive system \mca\ obtained from an additive system \mcb\ 
by the partition procedure described in 
Lemma~\ref{deBruijn:lemma:contraction} 
is called a \emph{contraction} of \mcb.  
(In~\cite{debr56},  de Bruijn called \mca\ a  \emph{degeneration}  
\index{degenerate additive system} of \mcb.) 
The set $I$ in Lemma~\ref{deBruijn:lemma:contraction} 
can be finite or infinite.  
If $I=J$ and $\sigma$ is a permutation of $J$ such that $J_i =\{\sigma(i) \}$
for all $i \in J$, then \mca\ and \mcb\ contain exactly the same sets.   
Thus, every additive system is a contraction of itself.  
An additive system \mca\ is a \emph{proper contraction} of \mcb\
if at least one set $A_i \in \mca$ is the sum of at least two sets in \mcb.

If $I = \{1\}$ and $J_1 = J$, then $A_1 = \N_0$. 
Thus, the additive system $(\N_0 )$ is a contraction of every additive system.  

The following Lemma shows that the contraction of a contraction is a contraction:

\bl     \label{deBruijn:lemma:ContractionTransitivity}
If $\mca = (A_i )_{i\in I}$,  $\mcb = (B_j )_{j \in J}$, 
and $\mcc  = (C_k)_{k \in K}$ are
additive systems such that \mca\ is a contraction of \mcb\ 
and \mcb\ is a contraction of \mcc, then \mca\ is a contraction of \mcc.  
\el

\begin{proof}
Because \mca\ is a contraction of \mcb, 
there exists a partition $\{J_i : i \in I\}$ of $J$ such that 
$A_i = \sum_{ j \in J_i} B_j$ for all $i \in I$. 
Because \mcb\ is a contraction of \mcc, 
there exists a partition $\{ K_j: j \in J \}$ of $K$ such that 
$B_j = \sum_{ k \in K_j} C_k$ for all $j \in J$.    
Then 
\[
A_i = \sum_{ j \in J_i} B_j  = \sum_{ j \in J_i} \sum_{k \in K_j} C_k 
=  \sum_{k \in L_i} C_k
\]
where 
\[
L_i = \bigcup_{j\in J_i} K_j \subseteq K
\]
and  
\[
\bigcup_{i \in I} L_i =   \bigcup_{i \in I} \bigcup_{j\in J_i} K_j 
= \bigcup_{j \in J} K_j = K.
\]
We shall show the sets in $\{L_i: i \in I \}$ are pairwise disjoint.

Let $k \in K$.  If $i_1, i_2 \in I$ and $k \in L_{i_1} \cap L_{i_2}$, 
then $k \in K_{j_1}$ for some $j_1 \in J_{i_1}$ 
and $k \in K_{j_2}$ for some $j_2 \in J_{i_2}$.  
Because  $\{K_j: j \in J \}$ is a set of pairwise disjoint sets 
and $K_{j_1} \cap K_{j_2} \neq \emptyset$, 
it follows that $j_1 = j_2$ and so $J_{i_1}\cap J_{i_2} \neq \emptyset$.
Because  $\{ J_i: i \in I \}$ is a set of pairwise disjoint sets, 
it follows that $i_1 = i_2$, and so the sets in $\{L_i: i \in I \}$ are pairwise disjoint. 
Thus, $\{L_i: i \in I \}$ is a partition of $K$, 
and the additive system  \mca\ is a contraction of \mcc.
This completes the  proof.
\end{proof}

Let \mca\ and \mcb\ be additive systems, let $r \in \N$, 
and let $(g_i)_{i\in [1,r]}$ be a finite sequence of integers $g_i \geq 2$.  
The expression ``\mca\ is a contraction of  \mcb\ dilated by 
 $(g_i)_{i\in [1,r]}$'' means that \mca\ is  the additive system
 obtained by first dilating \mcb\ by $(g_i)_{i\in [1,r]}$ and then contracting the dilated system.  
It was not hard to prove that a ``dilation of a dilation'' is a dilation 
(Lemma~\ref{deBruijn:lemma:DilationTransitivity}) or
that a ``contraction of a contraction is a contraction'' 
(Lemma~\ref{deBruijn:lemma:ContractionTransitivity}).
It is more challenging to prove that a ``contraction of a dilation of a contraction of a dilation'' is a contraction of a dilation.

 \bl    \label{deBruijn:lemma:DilateContract}
Let $\mca$, $\mcb$, and $\mcc$ 
be additive systems.  
If the additive system \mca\ is a contraction of the additive system \mcb\ 
dilated by the finite sequence $(g_i)_{i\in [1,r]}$,
and if \mcb\ is a contraction of the additive system \mcc\
dilated by the finite sequence $(g'_j)_{j \in [1,s]}$,
then \mca\ is a contraction of the additive system \mcc\ dilated by $(g_i)_{i\in [1,r+s]}$,
where $g_{r+j} = g'_j$ for $j \in [1,s]$.
\el

\begin{proof}
See Appendix~\ref{deBruijn:appendix}.
\end{proof}

\bl    \label{deBruijn:lemma:DilateContract2}
Let $(\mca_i)_{i\in [0,n]}$ be a sequence of additive systems and let 
$(g_i)_{i\in [1,n]}$ be a finite sequence of integers $g_i \geq 2$ such that 
$\mca_{i-1}$ is a contraction of $\mca_i$ dilated by $g_i$ for all $i \in [1,n]$.
Then $\mca_0$ is a contraction of $\mca_n$ dilated by 
the sequence $(g_i)_{i \in [1,n]}$.
\el

\begin{proof}
This follows from Lemma~\ref{deBruijn:lemma:DilateContract} by induction on $n$.  
\end{proof}

\section{British number systems}
In this section we describe certain additive systems that de Bruijn called
\index{British number system}   \emph{British number systems}.
A British number system 
is an additive system constructed from an infinite 
sequence of integers according to the algorithm in the following theorem.  

\bt    \label{deBruijn:theorem:BNS}
Let $(g_i )_{i \in \N}$ be an infinite sequence of integers 
such that $g_i \geq 2$ for all $i \geq 1$.  
Let $G_0 = 1$ and, for $i \in \N$, let  $G_i = \prod_{j=1}^i g_j$ and 
\[
A_i  = \{ 0, G_{i-1}, 2 G_{i-1},\ldots, (g_i-1)G_{i-1} \} = G_{i-1} \ast [0,g_i).
\]
Then $\mca = (A_i )_{i \in \N}$  is an additive system.  
\et

\begin{proof}
If $n \in \N_0$, then $n \in [0,G_r)$ for some sufficiently large integer $r$.  
By Theorem~\ref{deBruijn:theorem:DivisionAlgorithm}, 
there exist unique integers $a_i \in A_i$ for 
$i = 1,\ldots, r$ such that 
$n = \sum_{i=1}^r a_i \in \bigoplus_{i=1}^r A_i.$
Because $a \geq G_r$ for all 
$a \in \left( \bigcup_{i \in \N \setminus [1,r]} A_i \right)\setminus \{0\}$,
it follows that $n$ has a unique representation 
in the form $ \sum_{i \in \N} a_i $
with $a_i \in A_i$ for all $i \in \N$, and so 
$\mca = (A_i )_{i \in \N}$  is an additive system.  
\end{proof}

We write that the  sequence $(g_i )_{i \in \N}$ 
\emph{generates} the British 
number system \mca\ if \mca\ is constructed from $(g_i )_{i \in \N}$ 
according to the algorithm in Theorem~\ref{deBruijn:theorem:BNS}.

\bl
If $(g_i )_{i \in \N}$ generates the British number system $\mca = (A_i )_{i \in \N}$
and if $(g'_i )_{i \in \N}$  generates the British number system $\mca' = (A'_i )_{i \in \N}$, 
then $\mca = \mca'$ if and only if $g_i = g'_i$ for all $i \in \N$.  
\el

Thus, there is a one-to-one correspondence between British number systems 
and integer sequences $(g_i )_{i \in \N}$ 
satisfying $g_i \geq 2$ for all $i\in \N$.  

\begin{proof}
If $A_1 = A'_1$, then $[0,g_1) = [0,g'_1)$ and so $g_1 = g'_1$.  
If $r \geq 2$ and $g_i = g'_i$ for all $i  \leq r-1$, then 
\[
G_{r-1} = \prod_{i\in [1,r-1]}g_i = \prod_{i\in [1,r-1]}g'_i = G'_{r-1}.
\]
If $A_r = A'_r$, then  
\[
G_{r-1}\ast [0,g_r) = G'_{r-1}\ast [0,g'_r) = G_{r-1}\ast [0,g'_r) 
\]
and so $g_r = g'_r$.  
If $\mca = \mca'$, then it follows by induction that $g_i = g'_i$ for all $i \in \N$.
\end{proof}

de Bruijn's theorem is that 
\emph{every} additive system is a contraction of a British number system.  
The proof depends on the following fundamental lemma.

\bl     \label{deBruijn:lemma:deBruijn}
Let  $\mca = (A_i)_{i\in I}$ be an additive system.   
If $|I| \geq 2$, then there exist $i_1 \in I$, an integer $g \geq 2$,  
and a family of sets $\mcb = (B_i)_{i\in I}$ such that 
\[
A_{i_1} = [0, g)  \oplus g\ast B_{i_1}
\]
and, for all $i \in I \setminus \{i_1\}$,
\[
A_i = g\ast B_i.
\]
If $B_{i_1} = \{ 0\}$, then $\mcb = (B_i)_{i\in I\setminus \{i_1\}}$ is an additive system, 
and \mca\ is the dilation of the additive system \mcb\ by the integer $g$. 
If $B_{i_1} \neq \{ 0\}$, then $\mcb = (B_i)_{i\in I}$ is an additive system 
and \mca\ is a contraction of the additive system \mcb\ dilated by $g$.
\el

\begin{proof}
The inequality $|I| \geq 2$ implies that $A_i \neq \N_0$ for all $i\in I$.  
Because $1 \in \sum_{i\in I} A_i$, it follows that $1 \in A_{i_1}$ for some $i_1 \in I$.  
Because $A_{i_1} \neq \N_0$, there is a smallest positive integer $g$ 
such that $g \notin A_{i_1}$.   
Then $g \geq 2$ and $[0,g) \subseteq A_{i_1}$.
The sets in the family $(A_i \setminus \{ 0\})_{i\in I}$ are pairwise disjoint, 
and so $\{ 1,\ldots, g-1\} \cap A_i = \emptyset$ for all $i \in I\setminus \{ i_1\}$.  

We have $g = \sum_{i\in I} a_i \in \sum_{i\in I} A_i$, 
with $0 \leq a_i \leq g$ for all $i \in I$.
If  $1 \leq a_{i_1} \leq g-1$, then there must exist $j \in I\setminus \{i_1\}$ 
such that $1 \leq a_j \leq g-1$, which is absurd.  
Therefore, $a_{i_2}= g$ for some $i_2\in I\setminus \{i_1\}$ and 
$a_i = 0$ for all $i \in I\setminus \{ i_2\}$.  

Let $r \in \{1,2,\ldots, g-1\}$.  Then 
\[
r+g \in A_{i_1} + A_{i_2} \subseteq \sum_{i\in I} A_i.
\]
Because the representation of an integer in $\sum_{i\in I} A_i$ is unique, 
it follows that $r+g \notin A_i$ for all $i\in I.$

We shall prove that for every nonnegative integer $k$ the following holds:  
\benum
\item[(i)]
$[kg+1, (k+1)g)  \cap \bigcup_{i\in I\setminus \{ i_1 \} }A_i = \emptyset$,
\item[(ii)]
If  $[kg, (k+1)g) \cap A_{ i_1} \neq \emptyset$, then $ [kg, (k+1)g)  \subseteq A_{i_1}$.
\eenum
The proof is by induction on $k$ .  Statements~(i) and~(ii) already been verified 
for $k=0$ and $k = 1$.  
Let $k \geq 2$ and assume that statements~(i) and~(ii) are true 
for all nonnegative integers $k' < k$.  

For each $i \in I$ there exists $a_i \in A_i$ such that $0 \leq a_i \leq kg$ and 
\[
kg = \sum_{i \in I} a_i= a_{i_1} + \sum_{i \in I \setminus \{ i_1 \} } a_i.
\] 
By the induction hypothesis, $k' g + r \notin  \bigcup_{i\in I\setminus \{ i_1 \} }A_i$ 
for all $k' \in [0,k)$ and $r \in [1,g)$. 
Therefore, $a_i \equiv 0 \pmod{g}$ for all $i \in I \setminus \{ i_1 \}$,  and so $a_{i_1} \equiv 0 \pmod{g}$.

There are two cases.  
In the first case we have $kg \notin A_{i_1}$, and so  
$a_{i_1} = k' g$ for some nonnegative integer $k'< k$.  
By the induction hypothesis, $a_{i_1} + r= k' g + r \in A_{i_1}$  
for all $r \in [1,g)$, 
and so 
\[
kg + r = (a_{i_1} + r)  + \sum_{i \in I \setminus \{ i_1 \} } a_i \in \sum_{i \in I} A_i.
\]
Because the integer $kg+r$ has a unique representation  
in the sumset $\sum_{i \in I} A_i$, 
it follows that $kg+r \notin \bigcup_{i\in I} A_i$   for all $r \in [1,g)$.  

In the second case we have $kg \in A_{i_1}$.   
Because  $g \in A_{i_2}$, we have 
\[
(k+1)g = kg + g \in A_{i_1}+ A_{i_2}\subseteq \sum_{i\in I} A_i.
\] 
Let $r \in \{1,2,\ldots, g-1\} $. 
Because $\{1,2,\ldots, g-1\} \subseteq  A_{i_1}$, it follows that $g-r \in A_{i_1}$.  
If $kg+r \in A_{i_3}$ for some $r \in \{1,2,\ldots, g-1\} $ and $i_3 \neq i_1$, then
\[
(k+1)g = (g-r) + (kg+r) \in A_{i_1}+ A_{i_3} \subseteq \sum_{i\in I} A_i.
\]
This gives two distinct representations of $(k+1)g$ in $\sum_{i \in I} A_i$, which is absurd.  Therefore, $kg+r \notin A_i$ for all $i \in I \setminus \{ i_1 \}$.  
Thus, if $a_i \in A_i$ for $i \in I \setminus \{ i_1\}$ and $a _i < (k+1)g$, 
then $a_i \equiv 0 \pmod{g}.$   Writing 
\[
kg+r = a_{i_1} + \sum_{i\in I \setminus \{i_1\} } a_i \in \sum_{i\in I}A_i
\]
we conclude  that there exists a nonnegative integer $\ell \leq k$ such that  
\[
a_{i_1} = \ell g + r \in A_{i_1}
\]
and
\[
\sum_{i\in I \setminus \{i_1\} } a_i  = (k-\ell) g.
\]
If $\ell < k$, then the induction hypothesis implies that $\ell g \in A_{i_1}$ and so
\[
\ell g +  \sum_{i\in I \setminus \{i_1\} } a_i = kg,
\]
which is impossible since $kg \in A_{i_1}$.  
Therefore, $\ell = k$ and  $kg+r \in A_{i_1}$ 
for all $r \in [0, g)$  
This completes the induction.  

For each $i \in I$, let $B_{i} = \{k\in \N_0: kg \in A_{i} \}$.  Then 
\beq   \label{deBruijn:Ai1}
A_{i_1} = [0,g) \oplus g\ast B_{i_1}
\eeq
and, for every $i \in I \setminus \{ i_1 \}$, 
\[
A_i = g \ast B_i.
\]
Let $n \in N_0$.  There is a unique sequence of integers 
$(b_i)_{i\in I}$ with $b_i \in B_i$ for all $i \in I$ such that
\[
1+gn = (1+gb_{i_1}) + \sum_{i \in I \setminus \{i_1\}} gb_i \in \sum_{i\in I} A_1.
\]
It follows that  $n = \sum_{i \in I} b_i \in \sum_{i\in I}B_i$.     
If $B_{i_1} = \{0\}$, then $A_{i_1} = [0,g-1)$ 
and   $\mcb = (b_i)_{i \in I\setminus \{ i_1\}}$ is an additive system. 
Thus, \mca\ is the dilation of the additive system \mcb\ by the integer $g$.  

If $B_{i_1} \neq \{0\}$, then $\mcb = (b_i)_{i \in I}$ is an additive system 
and the decomposition~\eqref{deBruijn:Ai1} shows that 
\mca\ is a contraction of the additive system \mcb\ dilated by the integer $g$.
This completes the proof.
\end{proof}

We can now prove de Bruijn's theorem.

\bt   \label{deBruijn:theorem:deBruijn}
Every additive system is a British number system 
or a proper contraction of a British number system.
\et

\begin{proof}
Let $\mca = (A_i)_{i\in I}$ be an additive system, where, 
as usual, we assume that $I \cap \N = \emptyset$.  
If $|I| = 1$, then the additive system \mca\ consists of the single set $\N_0$, 
and $\N_0$ is a proper contraction of every British number system.

Let $\mca = \mca_0$.  
If $|I| \geq 2$, then Lemma~\ref{deBruijn:lemma:deBruijn} produces 
an additive system $\mca_1 = (A_{i,1})_{i\in I_1}$, with $I_1 \subseteq I$, 
and an integer $g_1 \geq 2$, such that $\mca_0$ is a contraction of $\mca_1$ 
dilated by $g_1$.  

Let $r \geq 1$, and suppose that we have constructed a sequence
$(g_i)_{i\in [1,r]}$ of integers $g_i \geq 2$ and a sequence of additive systems
$(\mca_i)_{i\in [0 ,r]}$ such that $\mca_{i-1}$ is a contraction of $\mca_i$ dilated by $g_i$
for all $i \in [1,r]$.  If $\mca_r = (A_{i,r})_{i\in I_r}$ and $|I_r| \geq 2$, 
then there is an additive system $\mca_{r+1} = (A_{i,1})_{i\in I_{r+1}}$, 
with $I_{r+1} \subseteq I_r \subseteq I$, and an integer $g_{r+1} \geq 2$ 
such that $\mca_r$ is a contraction of $\mca_{r+1}$ dilated by $g_{r+1}$.  

There are two cases.  
In the first case,  the process of constructing $\mca_{r+1}$ from $\mca_r$ 
terminates after $n$ steps.  This means that, after constructing the finite sequence 
of additive systems $(\mca_i)_{i\in [0,n]}$, 
we obtain $\mca_n = (A_{i,n})_{i\in I_n}$ with $|I_n| = 1$, 
that is, $\mca_n$ is the additive system consisting only of the set $\N_0$.
By Lemma~\ref{deBruijn:lemma:DilateContract2}, 
$\mca$ is a contraction of a dilation of $(\N_0)$ by the sequence $(g_i)_{i \in [1,n]}$. 
Because $(\N_0)$ is a contraction of every British number system, 
it follows that \mca\ is also a contraction of a British number system. 

In the second case, the process of constructing $\mca_{r+1}$ from $\mca_r$ 
never terminates, and we obtain an infinite sequence $(\mca_i)_{i\in \N}$
of additive systems and an infinite sequence $(g_i)_{i\in \N}$ of integers $g_i \geq 2$
such that $\mca_{i-1}$ is a contraction of the dilation of $\mca_i$ by $g_i$ 
for all $i \in \N$.  
By Lemma~\ref{deBruijn:lemma:DilateContract2}, we know that, 
for every positive integer $n$, the additive system $\mca$ 
is a contraction of $\mca_n$ dilated by the sequence $(g_i)_{i \in [1,n]}$. 

Recall that the additive system $\mca_n = (A_{i,n})_{i\in I_n}$ 
dilated by the sequence $(g_i)_{i \in [1,n]}$
consists of the sets $G_{i-1}\ast [0,g_i)$ for $i \in [1,n]$ 
and  $G_n \ast A_{i,n}$ for $i \in I$.

Let $\mca^{\flat}$ be the British number system generated 
by the infinite sequence $(g_i)_{i \in \N}$. 
We must prove that \mca\ is a contraction of $\mca^{\flat}$.  
Equivalently, we must construct a partition $(L_i)_{i\in I}$ of \N\ 
into pairwise disjoint nonempty sets such that 
\beq    \label{deBruijn:Li}
A_i = \sum_{n \in L_i} G_{n -1}\ast [0,g_{n})
\eeq
for all $i \in I$.  
Let 
\[
L_i = \{ n \in \N: G_{n-1} \in A_i \}
\]
and 
\[
I^{\flat} = \{i\in I: L_i \neq \emptyset \}.
\]

Let $n \in \N$.  
The additive system \mca\ is a contraction of the additive system $\mca_{n}$, 
and $G_{ n - 1}\ast [0,g_n)$ is a set in $\mca_n$.    
Therefore, the  set $G_{ n- 1}\ast [0,g_n)$ is  a summand 
in some set $A_i$ in \mca.  Because 
\[
G_{n -1} \in G_{ n - 1}\ast [0,g_n) \subseteq A_i
\]
it follows that $n \in L_i$ and so $\N = \bigcup_{i\in I^{\flat}} L_i$.
Because the sets $(A_i)_{i\in I}$ are pairwise disjoint, 
it follows that there is a unique $i \in I^{\flat}$ 
such that $G_{n-1} \in A_i$ and $n \in L_i$,
and so $(L_i)_{i\in I^{\flat}}$ is a partition of \N\ into nonempty, pairwise disjoint sets.

Let $i \in I$ and $x \in A_i\setminus \{ 0 \}$.  
Then $1 \leq x < G_N$ for some $N \in \N$.  
Because \mca\ is a contraction of $\mca_N$, 
the set $A_i$ is a sum of sets of the form $G_{n-1}\ast [0,g_n)$ with $n \in [1,N]$ 
and sets all of whose positive elements are greater than or equal to $G_N$.  
It follows that there is a nonempty subset $J$ of $[1,N]$ 
such that $x = \sum_{n \in J}G_{n-1}x_n$, 
with $x_n \in [1,g_n)$ and $G_{n-1} x_n \in A_i$ for all $n \in J$.  
This is possible only if $G_{n-1} \in G_{n-1}\ast [0,g_n) \subseteq A_i$ 
for all $n \in J$, and so $J \subseteq L_i$ and 
$x \in  \sum_{n \in L_i} G_{n -1}\ast [0,g_{n})$, that is,
\beq   \label{deBruijn:Li-1}
A_i \subseteq  \sum_{n \in L_i} G_{n -1}\ast [0,g_{n}).
\eeq
Moreover, $L_i \neq \emptyset$ implies that $i \in I^{\flat}$ and so $I^{\flat} = I$.

Conversely, if $i \in I$ and $ n \in L_i$, 
then $ G_{n - 1}\ast [0,g_n) \subseteq A_i $ and so 
\beq       \label{deBruijn:Li-2}
 \sum_{ n \in L_i} G_{ n -1}\ast [0,g_n) \subseteq A_i.
\eeq
The set inclusions~\eqref{deBruijn:Li-1} and~\eqref{deBruijn:Li-2}
imply~\eqref{deBruijn:Li}.  
This proves that \mca\ is a contraction of the  British number system $\mca^{\flat}$.
\end{proof}

\section{Remarks and open problems}
\textbf{Remark 1:}  
The set $A$ of integers is \emph{decomposable} 
\index{decomposable set} 
if there exist sets $B$ and $C$ such that $|B| \geq 2$, 
$|C| \geq 2$, and $A=B \oplus C$. 
An \emph{indecomposable set} is a set that does not decompose. 
An \emph{indecomposable additive system} is an  additive system 
$\mca = (A_i)_{i\in I}$ in which every set $A_i$ is indecomposable.  
Equivalently, an indecomposable additive system is an additive system 
that is not a proper contraction of another additive system.  
The following result  classifies indecomposable 
additive systems.

\bt[Nathanson~\cite{nath13a}]
Every infinite sequence of prime numbers generates an indecomposable 
British number system, and every indecomposable additive system is a British number system generated by an infinite sequence of prime numbers.  
There is a one-to-one correspondence between infinite sequences 
of prime numbers and indecomposable British number systems.  
Moreover, every additive system is either indecomposable or  a contraction of an  
indecomposable system.  
\et

\textbf{Remark 2:}  
Let $X$ be a nonempty set.  The \emph{free monoid}  on $X$ is the set $\mcm(X)$ 
consisting of all finite sequences of elements of $X$, 
and also an element $e$ (the ``empty sequence''), 
with the binary operation of \emph{concatenation}.  
We define the product of the nonempty sequences $(g_i)_{i\in [1,r]}$ 
and $(g'_j)_{j \in [1,s]}$ as follows:  
\[
(g_i)_{i\in [1,r]} \ast (g'_j)_{j \in [1,s]} = (g''_k)_{k \in [1,r+s]}
\]
where 
\[
g''_k = 
\begin{cases}
g_k & \text{ if } k \in [1,r] \\
g'_{k-r} & \text{ if } k \in [r+1,r+s]
\end{cases}
\]
and we define $e e = e$ and 
\[
(g_i)_{i\in [1,r]} e = e (g_i)_{i\in [1,r]}  = (g_i)_{i\in [1,r]}.
\]
The isomorphism class of the free monoid $\mcm(X)$ 
depends only on the cardinality of $X$.  
Lemma~\ref{deBruijn:lemma:DilationTransitivity} states that the free monoid 
on the set $\N\setminus \{1\}$ acts  by dilation on the set of additive systems.

\textbf{Remark 3:}  
Additive systems for the nonnegative integers are part of the  general 
study of sumsets.  If $A$ and $B$ are sets of integers, then their sumset is the set
$A+B = \{a+b:a\in A \text{ and } b \in B\}$.  
It is, in general, difficult to determine if a  set of integers is a sumset 
or ``almost'' a sumset, or to determine if a set is decomposable.  
Here are some open problems:  
Let $C$ be a nonempty finite or infinite set of integers.
\benum
\item
Do there exist sets $A$ and $B$ with $|A|\geq 2$, $|B|\geq 2$, and $A\oplus B=C$?
\item
Do there exist sets $A$ and $B$ with $|A|\geq 2$, $|B|\geq 2$, and  $A+ B=C$? 
\item
Do there exist sets $A$ and $B$ with $|A|\geq 2$ and $|B|\geq 2$ such that  
$A + B \subseteq C$ and $C \setminus (A+B)$ is ``small''? 
\item
Do there exist sets $A$ and $B$ with $|A|\geq 2$ and $|B|\geq 2$ such that  
$A + B \supseteq C$ and $(A+B) \setminus C$ is ``small''? 
\item
Does there exist a set $A$  with $|A|\geq 2$ and $A+A=C$? 
\item
Does there exist a set $A$  with $|A|\geq 2$ such that  
$A + A \subseteq C$ and $C \setminus (A+A)$ is ``small''? 
\item
Does there exist a set $A$  with $|A|\geq 2$ such that $A+A \supseteq C$ 
and $(A+A)\setminus C$ is ``small''? 
\eenum
These problems are related to Freiman's theorem~\cite{frei73} and 
other inverse problems in additive number theory 
(cf. Nathanson~\cite{nath96bb} and Tao and Vu~\cite{tao-vu}).

\textbf{Remark 4:}  
It is natural to investigate additive systems for the additive group \Z\ of integers, 
that is, sequences $(A_i)_{i\in I}$ of sets of integers such that $0 \in A_i$ 
and $|A_i| \geq 2$ for all $i \in I$, and 
$
\Z = \bigoplus_{i\in I} A_i.
$
For example, if $a_i = \varepsilon_i 2^{i-1}$ with $\varepsilon_i \in \{1, -1\}$ 
for all $i \in \N$, then $(\{0,a_i\})_{i \in \N}$ is an additive system for \Z\  
if and only if $\varepsilon_i = 1$ for infinitely many $i$ and $\varepsilon_i = -1$ 
for infinitely many $i$.  
The classification problem for additive systems for the integers is unsolved.  
Even the special case $A_i = \{0,a_i\}$ for all $i$ is difficult.  
de Bruijn~\cite{debr50} proved the following conjecture of T. Szele:  
If $(a_i)_{i \in \N}$  is an infinite sequence of nonzero integers such that 
$(\{0,a_i\})_{i \in \N}$ is an additive system for \Z, 
then there is a sequence $(d_i)_{i \in \N}$ of odd integers such that, 
after rearrangement, $a_i = 2^{i-1}d_i$ for all $i \in \N$.  

There are many interesting recent results about additive systems for \Z,  
for example, ~\cite{biro05,debr64,eige-haji08,eige-ito-pras04,ljuj12,schm-tull08,tijd95}.  
However, de Bruijn's remark at the end of his 1956 paper on $\N_0$ still accurately 
describes the current state of the problem:
``Some years ago the author~\cite{debr50} discussed various aspects 
of the analogous problem for number systems representing uniquely all 
integers (without restriction to nonnegative ones).  
That problem is much more difficult than the one dealt with above
[additive systems for $\N_0$], 
and it is still far from a complete solution.''

\textbf{Remark 5:}  
The interval identity $[0,mn) = [0,m) + m\ast [0,n)$, basic to the problem 
of additive systems for $\N_0$,  also led to the study 
of multiplication rules for quantum 
integers
(cf.~\cite{bori-nath-yang04,kont-nath06,nath03qa,nath11c,nguy10c}).

\textbf{Remark 6:}  
de Bruijn's paper~\cite{debr56} fills less than three pages.  
He uses but does not explicitly state  Lemma~\ref{deBruijn:lemma:DilateContract}, 
which is technically the most difficult step in the proof of the main result 
(Theorem~\ref{deBruijn:theorem:deBruijn}).  
After proving Lemma~\ref{deBruijn:lemma:deBruijn}, de Bruijn writes,   
``[Theorem~\ref{deBruijn:theorem:deBruijn}] 
easily follows by repeated application of the \ldots lemma.''
R. A. Rankin~\cite{rank} repeated this in his report on de Bruijn's paper in \emph{Mathematical Reviews}:  
``[Theorem~\ref{deBruijn:theorem:deBruijn}]  
follows from repeated applications of [the] lemma\ldots.''
Mathematicians, from the humblest  graduate student to the grandest Fields medalist, 
often don't bother to write out justifications for statements that are ``obvious''
or that ``easily follow'' from previously proved results.  But what is obvious to an author is not necessarily obvious to a reader (and sometimes the ``obvious'' is false).  
I prefer not to overindulge the virtue of brevity.

\appendix
\section{Proof of Lemma~\ref{deBruijn:lemma:DilateContract} }  \label{deBruijn:appendix}

\begin{proof}
For $k \in [1,s]$,  we define 
$
G'_k = \prod_{j =1}^k g'_j = \prod_{i=r+1}^{r+k} g_i 
$
and for $k \in [1,r+s]$,  we define 
$
G_k = \prod_{i=1}^k g_i.
$
If $k \in [1,s]$, then $G_rG'_k = G_{r+k}$ .
Let $G_0 = G'_0 =1$.

Let $\mcc'$ be the additive system $\mcc = (C_k)_{k\in K}$ 
dilated by $(g'_k)_{k \in [1,s]}$.  
We can assume that $K \cap \N = \emptyset$.
From the definition of dilation, we have  $\mcc'= (C'_k)_{k \in K'}$, where 
$
K' = [1,s]\cup K
$
and
\[
C'_k  = 
\begin{cases}
G'_{k - 1}\ast[0,g'_{k}) & \text{ if } k \in [1,s] \\
G'_{s} \ast C_k & \text{ if } k \in K.
\end{cases} 
\]
Let $\mcb = (B_j)_{j\in J}$ be a contraction of $\mcc'$,  
where $J \cap \N = \emptyset$.
This means that there is a partition 
$(K'_j)_{j\in J}$ of $K'$ such that $K'_j \neq \emptyset$  and 
\[
B_j = \bigoplus_{k\in K'_j} C'_k
\]
for all $j \in J$.

Let $\mcb'$ be the additive system $\mcb$ 
dilated by $(g_i)_{i\in [1,r]}$.  Then $\mcb'= (B'_j)_{j \in J'}$, where 
$
J'= [1,r] \cup J
$
and 
\begin{align*}
B'_j & = 
\begin{cases}
G_{j -1}\ast[0,g_j) & \text{ if } j \in [1,r] \\
G_r \ast B_j& \text{ if }  j \in J.
\end{cases} 
\end{align*}
Because $\mca = (A_i)_{i\in I}$ is a contraction of $\mcb'= (B'_j)_{j \in J'}$,  
there is a partition 
$(J'_i)_{i \in I}$ of $J'$ such that, for all $i \in I$, we have  $J'_i \neq \emptyset$ and 
\begin{align*}
A_i & = \bigoplus_{j \in J'_i} B'_j 
 = \left( \bigoplus_{j \in J'_i \cap [1,r]} G_{j -1}\ast [0,g_j ) \right)
\oplus  \left( \bigoplus_{j \in J'_i \setminus [1,r]} G_r \ast B_j  \right) \\
& = \left( \bigoplus_{j \in J'_i \cap [1,r]} G_{j -1}\ast [0,g_j ) \right)
\oplus  \left( \bigoplus_{j \in J'_i \setminus [1,r]} \bigoplus_{k\in K'_j} G_r \ast C'_k  \right).
\end{align*}
Note that $J'_i \setminus [1,r] \subseteq J$.  For $j \in J'_i \setminus [1,r]$, we have
\begin{align*}
\bigoplus_{k\in K'_j} G_r \ast C'_k 
& = \left( \bigoplus_{k\in K'_j \cap [1,s]} G_r \ast G'_{k - 1}\ast [0,g'_{k})  \right) 
\oplus  \left( \bigoplus_{k\in K'_j \setminus [1,s]} G_r \ast G'_s \ast C_k  \right) \\
& = \left( \bigoplus_{k\in K'_j \cap [1,s]} G_{r+k - 1}\ast [0,g_{r+k})  \right) 
\oplus  \left( \bigoplus_{k\in K'_j \setminus [1,s]} G_{r+s} \ast C_k  \right).
\end{align*}
It follows that
\begin{align*}
A_i  = 
& \left( \bigoplus_{j \in J'_i \cap [1,r]} G_{j -1}\ast [0,g_j ) \right)
\oplus  \left( \bigoplus_{j \in J'_i \setminus [1,r]} \bigoplus_{k\in K'_j \cap [1,s] } 
G_{r + k - 1}\ast [0,g_{r+k}) \right) \\
& \oplus  \left( \bigoplus_{j \in J'_i \setminus [1,r]} \bigoplus_{k\in K'_j \setminus [1,s]} G_{r+s} \ast C_k   \right).
\end{align*}
This is a decomposition of $A_i$ into a sum of sets.  We call these sets the 
\emph{summands} of $A_i$.  The summands of $A_i$ are pairwise distinct sets.

We must prove that $\mca = (A_i)_{i\in I}$  is a contraction of the additive system \mcc\ 
dilated by $(g_i)_{i=1}^{r+s}$.  
This dilated additive system can be written in the form 
$\mca^{\sharp} = (A^{\sharp}_k)_{k\in K^{\sharp} }$, 
where
$
K^{\sharp} = [1,r+s]  \cup K
$
and
\[
A^{\sharp}_k = 
\begin{cases} 
G_{k-1}\ast [0,g_k) & \text{ if } k \in [1,r+s] \\
G_{r+s}\ast C_k & \text{ if } k \in K.
\end{cases}
\]  
Every summand in $A_i$ is equal to $A^{\sharp}_k$ 
for some $k \in K^{\sharp}$.  
Thus, it suffices to show that for every $k \in K^{\sharp}$ there is a unique $i \in I$
such that $A^{\sharp}_k$ is a summand in $A_i$.  

The sets in the family $(J'_i)_{i\in I}$ partition $J' = [1,r]\cup J$.
Thus, for every $j\in [1,r]$ there is a unique $i \in I$ such that $j \in J'_i \cap [1,r]$,  
and so there is a unique $i \in I$ such that  $G_{j -1}\ast [0,g_j )$ is a summand in $A_i$.

Because the sets in the family $(K'_j)_{j\in J}$  partition $K' = [1,s]\cup K$, 
for every $k\in [1,s]$ there is a unique $j \in J$ such that $k \in K'_j \cap [1,s]$.  
The sets $\left( J'_i\setminus [1,r] \right)_{i\in I}$ partition $J$, and so 
there is a unique $i \in I$ such that $j \in J'_i\setminus [1,r]$.  
It follows that there is a unique $i \in I$ such that 
$G_{r + k - 1}\ast [0,g_{r+k})$ is a summand in $A_i$.  

Let $k \in K$.  There is a unique $j \in J$ such that $k \in K'_j \setminus [1,s]$,
and there is a unique $i \in I$ such that $j \in J'_i \setminus [1,r]$.  
It follows that there is a unique $i \in I$ such that $G_{r+s}\ast C_k$ 
is a summand in $A_i$.  
This proves that \mca\ is a contraction of the additive system $\mca^{\sharp}$.  
Indeed, defining 
\[
K^{\sharp}_i = \left( J'_i \cap [1,r] \right) \cup \left(  \bigcup_{j \in J'_i \setminus [1,r]} 
\left( r+ (K'_j \cap [1,s] ) \right) \right)
\cup \left(  \bigcup_{j \in J'_i \setminus [1,r]} K'_j \setminus [1,s] \right)
\]
we obtain a partition $\left( K^{\sharp}_i \right)_{i\in  [1,r+s] \cup I}$ of $K$ such that 
\[
A_i = \bigoplus_{k \in K^{\sharp}_i} A^{\sharp}_k
\]
for all $i \in I$.  
This completes the proof.
\end{proof}

\def\cprime{$'$} \def\cprime{$'$} \def\cprime{$'$} \def\cprime{$'$}
\providecommand{\bysame}{\leavevmode\hbox to3em{\hrulefill}\thinspace}
\providecommand{\MR}{\relax\ifhmode\unskip\space\fi MR }
\providecommand{\MRhref}[2]{%
  \href{http://www.ams.org/mathscinet-getitem?mr=#1}{#2}
}
\providecommand{\href}[2]{#2}

\end{document}